\documentclass[aop,MSNbibl,seceqn,citesort,dvips]{arximspdf}

\doi{10.1214/10-AOP591}
\volume{39}
\issue{4}
\pubyear{2011}
\firstpage{1502}
\lastpage{1527}

\makeatletter
\newcommand{\eqref}[1]{(\ref{#1})}
\newcommand{\Langle}{\langle\!\langle}
\newcommand{\Rangle}{\rangle\!\rangle}
\newcommand{\gfrac}[2]{#1/#2}
\newcommand{\vfrac}[2]{(#1/#2)}

\newproclaim{definition}{Definition}[section]

\newtheorem{theorem}{Theorem}[section]

\makeatother

\begin{document}
\begin{frontmatter}

\title{A probabilistic approach to Dirichlet problems of semilinear
elliptic PDE\textup{s} with singular coefficients}
\runtitle{A probabilistic approach to Dirichlet problems of PDE\textup{s}}

\begin{aug}
\author[A]{\fnms{Tusheng} \snm{Zhang}\corref{}\ead[label=e1]{tusheng.zhang@manchester.ac.uk}}
\runauthor{T. Zhang}
\affiliation{University of Manchester}
\address[A]{School of Mathematics\\
University of Manchester\\
Oxford Road\\
Manchester M13 9PL\\
United Kingdom\\
\printead{e1}} 
\end{aug}

\received{\smonth{10} \syear{2009}}
\revised{\smonth{5} \syear{2010}}

%
\begin{abstract}
In this paper, we prove that there exists a unique solution to the Dirichlet
boundary value problem for a general class of semilinear second order
elliptic partial
differential equations. Our approach is probabilistic. The theory of
Dirichlet processes
and backward stochastic differential equations play a crucial role.
\end{abstract}

%
\begin{keyword}[class=AMS]
\kwd[Primary ]{60H30}
\kwd[; secondary ]{35J25}
\kwd{31C25}.
\end{keyword}
\begin{keyword}
\kwd{Dirichlet processes}
\kwd{quadratic forms}
\kwd{Fukushima's decomposition}
\kwd{Dirichlet boundary value problems}
\kwd{backward stochastic differential equations}
\kwd{weak solutions}
\kwd{martingale representation theorem}.
\end{keyword}

\end{frontmatter}

\section{Introduction}\label{sec1}

In this paper, we will use probabilistic methods to solve the
Dirichlet boundary value problem for the semilinear
second order elliptic PDE of the following form:
%
\begin{equation}\label{0.1}
\cases{
\mathcal{A}u(x)=-f(x,u(x),\nabla u(x)),&\quad  $\forall x\in D$,\cr
u(x)|_{\partial D}=\varphi, &\quad $\forall
x\in\partial{D}$,
}
\end{equation}
where $D$ is a bounded domain in $R^d$.
The operator $\mathcal{A}$ is given by
%
\begin{equation}\label{0.0}
\qquad \mathcal{A}u={1\over2}\sum
_{i,j=1}^{d}{\partial\over{\partial x_i}}\biggl(a_{ij}(x)\,{\partial u\over
{\partial x_j}}\biggr)+
\sum_{i=1}^{d}b_i(x)\,{\partial u\over{\partial
x_i}}-\mbox{``$\operatorname{div}(\hat{b}u)$''}+q(x)u,
\end{equation}
where $a=(a_{i,j}(x))_{1\leq i,j\leq d}\dvtx D\to R^{d\times d}$ ($d>2$) is
a measurable, symmetric ma\-trix-valued function satisfying a uniform elliptic
condition, $b=(b_1,b_2,\ldots,b_d)$, $\hat{b}=(\hat{b}_1, \hat{b}_2,
\ldots,\hat{b}_d)\dvtx D\to R^d$ and $q\dvtx D\to R$ are merely
measurable functions belonging to some $L^p$ spaces, and $f(\cdot,
\cdot, \cdot)$ is a nonlinear function. The operator $\mathcal{A}$
is rigorously determined by the following quadratic form:
%
\begin{eqnarray}\label{0.00}
Q(u,v)&=&(-\mathcal{A}u, v)_{L^2(R^d)}\nonumber\\
&=&{1\over2}\sum_{i,j=1}^{d}\int_{R^d}a_{ij}(x)\,{\partial u\over
{\partial x_i}}{\partial v\over{\partial
x_j}}\,dx-\sum_{i=1}^{d}\int_{R^d}b_i(x)\,{\partial u\over{\partial
x_i}}v(x)\,dx\\
&&{}-\sum_{i=1}^{d}\int_{D}\hat{b}_i(x)u{\partial
v\over{\partial x_i}} \,dx-\int_{D}q(x)u(x)v(x) \,dx.\nonumber
\end{eqnarray}
We refer readers to \cite{GT,MR} and \cite{Tr} for details of the operator $\mathcal{A}$.

Probabilistic approaches to boundary value problems
of second order differential operators have been adopted by many
people. The earlier work went back as early as 1944 in \cite{Ka}.
See the books \cite{ChungZ,B} and references therein.
If $f=0$
(i.e., the linear case), and moreover $\hat{b}=0$, the solution $u$ to
problem (\ref{0.1}) can be solved by a Feynman--Kac formula
\[
u(x)= E_x \biggl[ \exp \biggl( \int_0^{\tau_D} q(X(s)) \,ds \biggr) \varphi(X(\tau_D)) \biggr]
\qquad \mbox{for } x\in D,
\]
where $X(t)$, $t\geq0$ is the diffusion process associated with the
infinitesimal generator
%
\begin{equation}\label{0.2}
L_1={1\over2}\sum_{i,j=1}^{d}{\partial\over{\partial
x_i}}\biggl(a_{ij}(x)\,{\partial\over{\partial
x_j}}\biggr)+\sum_{i=1}^{d}b_i(x)\,{\partial\over{\partial x_i}},
\end{equation}
$\tau_D$ is the first exit time of the diffusion process $X(t), t\geq
0$ from the domain $D$.
Very general results are obtained in the paper \cite{CZh} for this
case. When $\hat{b}\not=0$, ``$\operatorname{div}(\hat{b}\cdot)$'' in \eqref{0.0} is just a formal writing
because the divergence does not really exist for the merely
measurable vector field $\hat{b}$. It should be interpreted in the
distributional sense. It is exactly due to the nondifferentiability
of $\hat{b}$, all the previous known probabilistic
methods in solving the elliptic boundary value problems such as those in
\cite{CZh,B,Ka} and \cite{Ge} could not be applied. We
stress that the lower
order term $\operatorname{div}(\hat{b}\cdot)$
cannot be handled by Girsanov transform or Feynman--Kac transform either.
In a recent work \cite{CZ},
we show that the term $\hat{b}$ in fact can be tackled by the time-reversal
of Girsanov transform from the first exit time $\tau_D$ from $D$ by
the symmetric
diffusion $X^0$ associated with $L_0={1\over2}\sum
_{i,j=1}^{d}{\partial\over{\partial x_i}}(a_{ij}(x)\,{\partial\over
{\partial
x_j}})$, the symmetric part of $\mathcal{A}$. The
solution to equation (\ref{0.1}) (when $f=0$ ) is given by
%
\begin{eqnarray} \label{0.3}
u(x)&=& E_x^0 \biggl[\varphi(X^0(\tau_D)) \exp
\biggl\{\int_0^{\tau_D}\langle(a^{-1}b) (X^0(s)), dM^0(s)\rangle\nonumber\\
&&\hphantom{E_x^0 \biggl[\varphi(X^0(\tau_D)) \exp\biggl\{}
{}+ \biggl(\int_0^{\tau_D}\langle(a^{-1}\hat{b}) (X^0(s)), dM^0(s)\rangle
\biggr)\circ
r_{\tau_D}\nonumber\\[-8pt]\\[-8pt]
&&\hphantom{E_x^0 \biggl[\varphi(X^0(\tau_D)) \exp\biggl\{}
{}-\frac{1}{2}\int_0^{\tau_D}(b-\hat{b})a^{-1}(b-\hat{b})^{\ast}(X^0(s))\,ds\nonumber \\
&&\hspace*{186pt}
{}+\int_0^{\tau_D}q(X^0(s))\,ds \biggr\} \biggr],\nonumber
\end{eqnarray}
where $M^0(s)$ is the martingale part of the diffusion $X^0$, $r_t$
denotes the reverse operator, and $\langle\cdot, \cdot\rangle$ stands for the inner
product in $R^d$.

Nonlinear elliptic PDEs [i.e., $f\not=0$ in (\ref{0.1})] are generally very hard to solve. One can not
expect explicit expressions for the solutions. However, in recent
years backward stochastic differential equations (BSDEs) have been
used effectively to solve certain nonlinear PDEs. The general
approach is to represent the solution of the nonlinear equation
(\ref{0.1}) as the solution of certain BSDEs associated with the
diffusion process generated by the linear operator $\mathcal{A}$.
But so far, only the cases where $\hat{b}=0$ and $b$ being bounded
were considered. The main difficulty for treating the general
operator $\mathcal{A}$ in (\ref{0.0}) with $\hat{b}\not=0$, $q\not=0$ is that there are no associated diffusion processes anymore. The
mentioned methods used so far in the literature ceased to work.
Our approach is to transform the problem (\ref{0.1}) to a similar
problem for which the operator $\mathcal{A}$ does not have the
``bad'' term $\hat{b}$. See below for detailed description.

There exist many papers on BSDEs and their applications to nonlinear
PDEs. We mention some related earlier results. The first result on
probabilistic interpretation for solutions of semilinear parabolic
PDE's was obtained by Peng in \cite{P} and subsequently in \cite{PP2}. In \cite{DP}, Darling and
Pardoux obtained a viscosity solution to the Dirichlet problem for a
class of semilinear elliptic PDEs (through BSDEs with random
terminal time) for which the linear operator $\mathcal{A}$ is of the
form
\[
\mathcal{A}={1\over
2}\sum_{i,j=1}^{d}a_{ij}(x)\,{\partial^2\over{\partial
x_j\partial x_i}}+ \sum_{i=1}^{d}b_i(x)\,{\partial\over
{\partial x_i}},
\]
where $a_{ij}\in C_b^2(D)$ and $b\in C_b^1(D)$. BSDEs associated
with Dirichlet processes and weak solutions of semi-linear parabolic
PDEs were considered by Lejay in \cite{L} where the linear operator
$\mathcal{A}$ is assumed to be
\[
\mathcal{A}={1\over
2}\sum_{i,j=1}^{d}{\partial\over{\partial
x_i}}\biggl(a_{ij}(x)\,{\partial\over{\partial x_j}}\biggr)+
\sum_{i=1}^{d}b_i(x)\,{\partial\over{\partial x_i}},
\]
for bounded coefficients $a$ and $b$. BSDEs associated with
symmetric Markov processes and weak solutions of semi-linear
parabolic PDEs were studied by Bally, Pardoux and Stoica in
\cite{BPS} where the linear operator $\mathcal{A}$ is assumed to be
symmetric with respect to some measure $m$. BSDEs and solutions of
semi-linear parabolic PDEs were also considered by Rozkosz in
\cite{R} for the linear operator $\mathcal{A}$ of the form
\[
\mathcal{A}={1\over
2}\sum_{i,j=1}^{d}{\partial\over{\partial
x_i}}\biggl(a_{ij}(t,x){\partial\over{\partial x_j}}\biggr).
\]

Now we describe the contents of this paper in more details. Our
strategy is to transform the problem (\ref{0.1}) by a kind of
$h$-transform to a problem of a similar kind, but with an operator
$\mathcal{A}$ that does not have the ``bad'' term $\hat{b}$. The
first step will be to solve (\ref{0.1}) assuming $\hat{b}=0$. In
Section~\ref{sec2}, we introduce the Feller diffusion process $(\Omega, \mathcal{
F}, \mathcal{ F}_t, X(t), P_x, x\in R^d)$
whose infinitesimal generator is given by
%
\begin{equation}\label{0.4}
L_1={1\over2}\sum_{i,j=1}^{d}{\partial\over{\partial
x_i}}\biggl(a_{ij}(x)\,{\partial\over{\partial
x_j}}\biggr)+\sum_{i=1}^{d}b_i(x)\,{\partial\over{\partial x_i}}.
\end{equation}
In general, $X(t)$, $t\geq0$ is not a semimartingale. But it has a
nice martingale part $M(t)$, $t\geq0$. In this section, we prove a
martingale representation theorem for the martingale part $M(t)$,
which is crucial for the
study of BSDEs in subsequent sections. In Section~\ref{sec3}, we solve a
class of BSDEs associated with the martingale part $M(t)$, $t\geq0$:
%
\begin{equation}\label{0.5}
Y(t)=\xi+\int_t^Tf(s, Y(s), Z(s))\,ds-\int_t^T\langle Z(s),dM(s)\rangle .
\end{equation}
The random coefficient $f(t,y,z, \omega)$ satisfies a certain
monotonicity condition which is particularly fulfilled in the situation
we are interested. The
BSDEs with deterministic terminal time were solved first and then
the BSDEs with random terminal time were studied. In Section~\ref{sec4}, we
consider the Dirichelt problem for the second order differential
operator
%
\begin{equation}\label{0.6}
L_2={1\over2}\sum_{i,j=1}^{d}{\partial\over{\partial
x_i}}\biggl(a_{ij}(x)\,{\partial\over{\partial
x_j}}\biggr)+\sum_{i=1}^{d}b_i(x)\,{\partial\over{\partial
x_i}}+q(x),
\end{equation}
where $b_i\in L^p$ for some $p>d$ and $q\in L^{\beta}$ for some $\beta
>\frac{d}{2}$. We first solve the linear problem with a given function $F$
%
\begin{equation}\label{0.7}
\cases{ L_2 u= F ,&\quad in $D$, \cr
u= \varphi,&\quad on $ \partial D$,
}
\end{equation}
and then the nonlinear problem
%
\begin{equation}\label{0.8}
\cases{
L_2 u= -g(x, u(x), \nabla u(x)), &\quad in $D$, \cr
u= \varphi,&\quad on $\partial D$,
}
\end{equation}
with the help of BSDEs. Finally, in Section~\ref{sec5}, we study the Dirichlet
problem
%
\begin{equation}\label{0.9}
\cases{
\mathcal{A}u(x)=-f(x,u(x)),&\quad  $\forall x\in D$,\cr
u(x)|_{\partial D}=\varphi, &\quad $\forall
x\in\partial{D}$,
}
\end{equation}
where $\mathcal{A}$ is a general second order differential operator
given in (\ref{0.0}). We apply a
transform we introduced in \cite{CZ} to transform the above problem
to a problem like (\ref{0.8}) and then a reverse transformation will
solve the final problem.

\section{Preliminaries}\label{sec2}
 Let $\mathcal{A}$ be an elliptic operator
of the following general form:
\[
\mathcal{A}={1\over2}\sum_{i,j=1}^{d}{\partial\over{\partial
x_i}}\biggl(a_{ij}(x)\,{\partial\over{\partial x_j}}\biggr)+
\sum_{i=1}^{d}b_i(x)\,{\partial\over{\partial
x_i}}-\mbox{``$\operatorname{div}(\hat{b}\cdot)$''}+q(x),
\]
where $a=(a_{ij}(x))\dvtx D\to R^{d\times d}$ ($d>2$) is a measurable,
symmetric matrix-valued function which satisfies the uniform elliptic
condition
%
\begin{equation}\label{1.0}
\lambda|\xi|^2\leq\sum_{i,j=1}^{d}a_{ij}(x)\xi_i\xi_j\leq\Lambda
|\xi|^2,\qquad  \forall \xi\in R^d \mbox{ and } x\in D
\end{equation}
for some constant $\lambda, \Lambda>0$,
$b=(b_1,\ldots,b_d), \hat{b}=(\hat{b}_1,\ldots,\hat{b}_d)\dvtx D\to
R^d$ and $q\dvtx\break D\to R$ are measurable functions which could be
singular and such that
\[
|b|^2\in L^p(D),\qquad  |\hat{b}|^2\in L^p(D)\quad  \mbox{and}\quad  q\in L^p(D),
\]
for some $p>\frac{d}{2}$. Here $D$ is a bounded domain in
$R^d$ whose boundary is regular, that is, for every $x\in\partial D$,
$P(\tau_D^x=0)=1$, where
$\tau_D^x$ is the first exit time of a standard Brownian motion
started at $x$ from the domain $D$. Let $f\dvtx R^d\times R\times
R^d\rightarrow R$ be a measurable
nonlinear function. Consider the following nonlinear Dirichlet
boundary value problem:
%
\begin{equation}\label{1.01}
\cases{
\mathcal{A}u(x)=-f(x,u(x),\nabla u(x)),&\quad $ \forall x\in D$,\cr
u(x)|_{\partial D}=\varphi, &\quad $\forall
x\in\partial{D}$.
}
\end{equation}
Let $W^{1,2}(D)$ denote the usual Sobolev space of order one:
\[
W^{1,2}(D)=\{u\in L^2(D)\dvtx \nabla u\in L^2(D; R^d)\}.
\]

\begin{definition} We say that $u\in W^{1,2}(D)$ is a continuous, weak
solution of (\ref{1.01}) if:
\begin{itemize}[(iii)]
\item[(i)] for any $\phi\in W^{1,2}_{0}(D)$,
\begin{eqnarray*}
&&{1\over2}\sum_{i,j=1}^{d}\int_{D}a_{ij}(x)\,{\partial u\over{\partial
x_i}}\,{\partial\phi\over{\partial x_j}}\,dx-
\sum_{i=1}^{d}\int_{D}b_i(x)\,{\partial u\over{\partial
x_i}}\,\phi \,dx\\
&&\qquad {} -\sum_{i=1}^{d}\int_{D}\hat{b}_i(x)u\,{\partial
\phi\over{\partial x_i}}\, dx
-\int_{D}q(x)u(x)\phi \,dx=\int_{D}f(x,u,\nabla u )\phi\, dx,
\end{eqnarray*}

\item[(ii)] $u\in C(\bar{D})$,

\item[(iii)] $\lim_{y\rightarrow x}u(y)=\varphi
(x)$, $\forall x\in\partial{D}$.
\end{itemize}
\end{definition}

Next we introduce two diffusion processes which will be used later.\vadjust{\goodbreak}

Let $(\Omega, \mathcal{ F}, \mathcal{ F}_t, X(t), P_x, x\in R^d)$ be
the Feller diffusion process
whose infinitesimal generator is given by
%
\begin{equation}\label{1.1}
L_1={1\over2}\sum_{i,j=1}^{d}{\partial\over{\partial
x_i}}\biggl(a_{ij}(x)\,{\partial\over{\partial
x_j}}\biggr)+\sum_{i=1}^{d}b_i(x)\,{\partial\over{\partial x_i}},
\end{equation}
where $\mathcal{ F}_t$ is the completed, minimal admissible filtration
generated by $X(s)$, $s\geq0$. The
associated nonsymmetric, semi-Dirichlet form with $L_1$ is defined
by
%
\begin{eqnarray}\label{1.2}
\qquad Q_1(u,v)&=&(-L_1u, v)_{L^2(R^d)}\nonumber\\[-8pt]\\[-8pt]
&=&{1\over2}\sum_{i,j=1}^{d}\int_{R^d}a_{ij}(x)\,{\partial u\over
{\partial x_i}}\,{\partial v\over{\partial
x_j}}\,dx-\sum_{i=1}^{d}\int_{R^d}b_i(x)\,{\partial u\over{\partial
x_i}}\,v(x)\,dx.\nonumber
\end{eqnarray}
The process $X(t)$, $t\geq0$ is not a semimartingale in general.
However,
it is known (see, e.g., \cite{CZh,Fu,FOT} and \cite{LZ}) that the following
Fukushima's decomposition holds:
%
\begin{equation}\label{1.3}
X(t)=x+M(t)+N(t)\qquad  P_x\mbox{-a.s.},
\end{equation}
where $M(t)$ is a continuous square integrable martingale with sharp
bracket being given by
%
\begin{equation}\label{1.4}
\Langle M^i,M^j\Rangle_t=\int_0^ta_{i,j}(X(s))\,ds,
\end{equation}
and $N(t)$ is a continuous process of zero quadratic variation.
Later we also write $X_x(t)$, $M_x(t)$ to emphasize the dependence on
the initial value $x$. Let $\mathcal{ M}$ denote the space of square
integrable martingales w.r.t. the filtration $\mathcal{ F}_t$, $t\geq0$.
The following result is a martingale representation theorem whose proof is
a modification of that of Theorem A.3.20 in \cite{FOT}. It will play
an important role in our study of the backward stochastic
differential equations associated with the martingale part $M$.
\begin{theorem}
For any $L\in\mathcal{ M}$, there exist predictable processes $H_i(t),
i=1,\ldots,d$ such that
%
\begin{equation}\label{1.5}
L_t=\sum_{i=1}^d\int_0^t H_i(s)\,dM^i(s).
\end{equation}
\end{theorem}

\begin{pf}
 It is sufficient to prove
(\ref{1.5}) for $0\leq t\leq T$, where $T$ is an arbitrary, but
fixed constant $T$. Recall that $\mathcal{ M}$ is a Hilbert space w.r.t.
the inner product $(K_1, K_2)_\mathcal{ M}=E[\Langle K_1, K_2\Rangle _T]$, where
$\Langle K_1, K_2\Rangle $
denotes the sharp bracket of $K_1$ and $K_2$. Let $\mathcal{ M}_1$
denote the subspace of square integrable martingales of the form
(\ref{1.5}). Let $R_{\alpha}, \alpha>0$ be the resolvent operators
of the diffusion process\vadjust{\goodbreak} $X(t)$, \mbox{$t\geq0$}. Fix any $g\in C_b(R^d)$,
we know that $R_{\alpha}g\in D(L_1)$ and $L_1R_{\alpha}g=\alpha
R_{\alpha}g-g$. Moreover, it follows from \cite{FOT} and \cite{LZ}
that
\begin{eqnarray*}
R_{\alpha}g(X(t))-R_{\alpha}g(X(0))&=&\int_0^t\langle\nabla
R_{\alpha}g(X(s)), dM(s)\rangle\\
&&{}+\int_0^t(\alpha R_{\alpha}g-g)(X(s))\,ds.
\end{eqnarray*}
Hence,
\begin{eqnarray*}
J_t&:=&\int_0^t e^{-\alpha s}\langle\nabla
R_{\alpha}g(X(s)), dM(s)\rangle
\\
&\hspace*{3pt}=&e^{-\alpha t}R_{\alpha}g(X(t))-R_{\alpha}g(X(0))+\int_0^te^{-\alpha
s}g(X(s))\,ds
\end{eqnarray*}
is a bounded martingale that belongs to $\mathcal{ M}_1$.
The theorem will be proved if we can show that $\mathcal{
M}_1^{\perp}=\{0\}$. Take $K\in\mathcal{ M}_1^{\perp}$. Since $\mathcal{
M}_1$ is stable under stopping, by Lemma 2 in Chapter IV in
\cite{PEP}, we deduce $\Langle K, L\Rangle =0$ for all $L\in\mathcal{ M}_1$. In
particular, $\Langle K, J\Rangle =0$. From here, we can follow the same proof of
Theorem~A.3.20 in \cite{FOT}
to conclude $K=0$.
\end{pf}

We will
denote by $(\Omega, \mathcal{ F}^0, \mathcal{ F}_t^0, X^0(t), P^0_x,
x\in R^d)$ the diffusion process generated by
%
\begin{equation}\label{1.7}
L_0={1\over2}\sum_{i,j=1}^{d}{\partial\over{\partial
x_i}}\biggl(a_{ij}(x)\,{\partial\over{\partial
x_j}}\biggr).
\end{equation}
The corresponding Fukushima's decomposition is written as
\[
X^0(t)=x+M^0(t)+N^0(t),\qquad  t\geq0.
\]
For $v\in W^{1,2}(R^d)$, the Fukushima's decomposition for the
Dirichlet process $v(X^0(t))$
reads as
%
\begin{equation}\label{1.8}
v(X^0(t))=v(X^0(0))+M^v(t)+N^v(t),
\end{equation}
where $M^v(t)=\int_0^t\nabla v(X^0(s))\cdot dM^0(s)$, $N^v(t)$ is a
continuous process of zero energy (the zero energy part). See \cite{CFKZ1,CFKZ2,FOT} for details of symmetric Markov processes.

\section{Backward SDEs with singular coefficients}\label{sec3}
 Let $(\Omega, \mathcal{ F}, \mathcal{ F}_t)$
be the
probability space carrying the diffusion process $X(t)$ described in
Section~\ref{sec2}. Recall $M(t)$, $t\geq0$ is the martingale part of $X$. In
this section, we will study backward stochastic differential
equations (BSDEs) with singular coefficients associated with the
martingale part $M(t)$.

\subsection{BSDEs with deterministic terminal times}\label{sec3.1}

Let $f(s,y,z, \omega)\dvtx [0, T]\times R\times R^d\times\Omega
\rightarrow R$ be
a given progressively
measurable function. For simplicity, we omit the random parameter
$\omega$.
Assume that $f$ is continuous in $y$ and
satisfies:\vadjust{\goodbreak}
\begin{itemize}[(A.1)]
\item[(A.1)] $(y_1-y_2)(f(s,y_1,z)-f(s,y_2,z))\leq
-d_1(s)|y_1-y_2|^2$,
\item[(A.2)] $|f(s,y,z_1)-f(s,y,z_2)|\leq d_2 |z_1-z_2|$,
\item[(A.3)] $|f(s,y,z)|\leq|f(s,0,z)|+ K(s)(1+ |y|)$,
\end{itemize}
where $d_1(\cdot)$, $K(s)$ are a progressively
measurable stochastic process and $d_2$ is a constant.
Let $\xi\in L^2(\Omega, \mathcal{ F}_T, P)$. Let $\lambda$ be the constant
defined in (\ref{1.0}).

\begin{theorem}\label{thm3.1}
Assume $E [ e^{-\int_0^T2d_1(s)\,ds}|\xi|^2 ]<\infty$, $E[\int_0^T K(s)\,ds]<\infty$ and
\[
E \biggl[ \int_0^Te^{-\int_0^s2d_1(u)\,du}|f(s,0,0)|^2\,ds \biggr]<\infty.
\]
Then, there exists a unique ($\mathcal{ F}_t$-adapted) solution $(Y, Z)$
to the following BSDE:
%
\begin{equation}\label{2.1}
Y(t)=\xi+\int_t^Tf(s, Y(s), Z(s))\,ds-\int_t^T\langle Z(s),dM(s)\rangle ,
\end{equation}
where $Z(s)=(Z_1(s),\ldots,Z_d(s))$.
\end{theorem}

\begin{pf}
We first prove the uniqueness.
Set $d(s)=-2d_1(s)$. Suppose $(Y^1(t), Z^1(t))$ and $(Y^2(t),
Z^2(t))$ are two solutions to equation (\ref{2.1}). Then
%
\begin{eqnarray}\label{2.2}
&&d(|Y^1(t)-Y^2(t)|^2)\nonumber\\
&&\qquad =-2\bigl(Y^1(t)-Y^2(t)\bigr)\bigl(f(t,Y^1(t),
Z^1(t))-f(t,Y^2(t), Z^2(t))\bigr)\,dt\nonumber\\[-8pt]\\[-8pt]
&&\qquad\quad  {}+2\bigl(Y^1(t)-Y^2(t)\bigr)\langle Z^1(t)-Z^2(t),dM(t)\rangle \nonumber\\
&&\qquad\quad  {}+\bigl\langle a(X(t))\bigl(Z^1(t)-Z^2(t)\bigr), Z^1(t)-Z^2(t)\bigr\rangle \,dt.\nonumber
\end{eqnarray}
By the chain rule, using the assumptions (A.1), (A.2) and Young's
inequality, we get
\begin{eqnarray}\label{2.3}
&&|Y^1(t)-Y^2(t)|^2e^{\int_0^td(s)\,ds}\nonumber\\
&&\quad {}+\int_t^Te^{\int
_0^sd(u)\,du}\langle a(X(s))\bigl(Z^1(s)-Z^2(s)\bigr),
Z^1(s)-Z^2(s)\rangle\,ds \nonumber\\
&&\qquad =-\int_t^T e^{\int_0^sd(u)\,du}|Y^1(s)-Y^2(s)|^2d(s)\,ds\nonumber\\
&&\qquad \quad {} +2\int_t^T e^{\int_0^sd(u)\,du}\bigl(Y^1(s)-Y^2(s)\bigr)\nonumber \\
&&\qquad \quad \hphantom{{} +2\int_t^T }
{}\times\bigl(f(s,Y^1(s),
Z^1(s))-f(s,Y^2(s), Z^2(s))\bigr)\,ds\nonumber\\
&&\qquad \quad {}-2\int_t^T e^{\int_0^sd(u)\,du}\bigl(Y^1(s)-Y^2(s)\bigr) \langle Z^1(s)-Z^2(s),
dM(s) \rangle\nonumber
\\
&&\qquad \leq-\int_t^T e^{\int_0^sd(u)\,du}|Y^1(s)-Y^2(s)|^2d(s)\,ds\\
&&\qquad\quad  {} -2\int_t^T
e^{\int_0^sd(u)\,du}|Y^1(s)-Y^2(s)|^2d_1(s)\,ds\nonumber\\
&&\qquad\quad  {} +2\int_t^T
e^{\int_0^sd(u)\,du}d_2|Y^1(s)-Y^2(s)||Z^1(s)-Z^2(s)|\,ds\nonumber\\
&&\qquad\quad  {} -2\int_t^T
e^{\int_0^sd(u)\,du}\bigl(Y^1(s)-Y^2(s)\bigr) \langle Z^1(s)-Z^2(s),dM(s) \rangle\nonumber\\
&&\qquad \leq C_{\lambda}\int_t^T
e^{\int_0^sd(u)\,du}|Y^1(s)-Y^2(s)|^2\,ds\nonumber\\
&&\qquad\quad  {} +\frac{1}{2}\int_t^T
e^{\int_0^sd(u)\,du} \bigl\langle a(X(s))\bigl(Z^1(s)-Z^2(s)\bigr), \bigl(Z^1(s)-Z^2(s)\bigr) \bigr\rangle\,ds\nonumber
\\
&&\qquad\quad  {} -2\int_t^T
e^{\int_0^sd(u)\,du}\bigl(Y^1(s)-Y^2(s)\bigr) \langle Z^1(s)-Z^2(s),dM(s) \rangle.\nonumber
\end{eqnarray}
Take expectation in above inequality to get
\[
E \bigl[|Y^1(t)-Y^2(t)|^2e^{\int_0^td(s)\,ds} \bigr]
 \leq C_{\lambda}\int_t^TE\bigl[
e^{\int_0^sd(u)\,du}|Y^1(s)-Y^2(s)|^2\bigr]\,ds.
\]
By Gronwall's inequality,
we conclude $Y^1(t)=Y^2(t)$ and hence $Z^1(t)=Z^2(t)$ by (\ref{2.3}).

Next, we prove the existence. Take an even,
nonnegative function $\phi\in C_0^{\infty}(R)$ with
$\int_R\phi(x)\,dx=1$. Define
\[
f_n(t,y, z)=\int_Rf(t,x,z)\phi_n(y-x)\,dx,
\]
where $\phi_n(x)=n\phi(nx)$. Since $f$ is continuous in $y$, it
follows that $f_n(t,y, z)\rightarrow f(t,y, z)$ as $n\rightarrow
\infty$. Furthermore, it is easy to see that for every $n\geq1$,
%
\begin{equation}\label{2.4}
\vert f_n(t,y_1, z)-f_n(t,y_2, z)\vert\leq C_n |y_1-y_2|,\qquad  y_1, y_2 \in R,
\end{equation}
for some constant $C_n$.
Consider the following BSDE:
%
\begin{equation}\label{2.5}
Y_n(t)=\xi+\int_t^Tf_n(s, Y_n(s), Z_n(s))\,ds-\int_t^T \langle Z_n(s),
dM(s) \rangle.
\end{equation}
In view of (\ref{2.4}) and the assumptions (A.2), (A.3), it is known
(e.g., \cite{PP1}) that the above equation admits a unique solution $(Y_n,
Z_n)$. Our aim now is to show that there exists a convergent
subsequence $(Y_{n_k}, Z_{n_k})$. To this end, we need some
estimates. Applying It\^o's formula, in view of assumptions
(A.1)--(A.3) it follows that
%
\begin{eqnarray}\label{2.6}
&&|Y_n(t)|^2e^{\int_0^td(s)\,ds}+\int_t^Te^{\int_0^sd(u)\,du}\langle a(X(s))Z_n(s),
Z_n(s)\rangle\,ds \nonumber\\
&&\qquad =|\xi|^2e^{\int_0^Td(s)\,ds}-\int_t^T e^{\int
_0^sd(u)\,du}Y_n^2(s)d(s)\,ds\nonumber\\
&&\qquad \quad {} +2\int_t^T e^{\int_0^sd(u)\,du}Y_n(s)f_n(s,Y_n(s), Z_n(s))\,ds\nonumber\\
&&\qquad \quad {}-2\int_t^T e^{\int_0^sd(u)\,du}Y_n(s)\langle Z_n(s),dM(s)\rangle\nonumber\\
&&\qquad \leq|\xi|^2e^{\int_0^Td(s)\,ds}-\int_t^T e^{\int
_0^sd(u)\,du}Y_n^2(s)d(s)\,ds\nonumber\\
&&\qquad \quad {} -2\int_t^T e^{\int_0^sd(u)\,du}d_1(s)Y_n^2(s)\,ds+2C\int_t^T e^{\int
_0^sd(u)\,du}|Y_n(s)||Z_n(s)|\,ds\nonumber\\
&&\qquad \quad {}+2\int_t^T
e^{\int_0^sd(u)\,du}|Y_n(s)|f(s,0,0)\,ds\nonumber\\
&&\qquad \quad {}-2\int_t^T e^{\int_0^sd(u)\,du}Y_n(s)\langle Z_n(s),dM(s)\rangle\\
&&\qquad \leq|\xi|^2e^{\int_0^Td(s)\,ds}+C_{\lambda}\int_t^T e^{\int
_0^sd(u)\,du}Y_n^2(s)\,ds\nonumber\\
&&\qquad \quad {}+\frac{1}{2}\int_t^T e^{\int_0^sd(u)\,du}\langle a(X(s))Z_n(s), Z_n(s)\rangle\,ds
\nonumber\\
&&\qquad \quad {}+\int_t^T e^{\int_0^sd(u)\,du}Y_n^2(s)\,ds+\int_t^T
e^{\int_0^sd(u)\,du}|f(s,0,0)|^2\,ds\nonumber\\
&&\qquad \quad {}-2\int_t^T e^{\int_0^sd(u)\,du}Y_n(s)\langle Z_n(s),dM(s)\rangle.\nonumber
\end{eqnarray}
Take expectation in (\ref{2.6}) to obtain
%
\begin{eqnarray}\label{2.7}
&&E\bigl[|Y_n(t)|^2e^{\int_0^td(s)\,ds}\bigr]+\frac{1}{2}E\biggl[\int_t^Te^{\int
_0^sd(u)\,du}\langle a(X(s))Z_n(s),
Z_n(s)\rangle\,ds \biggr]\nonumber\\
&&\qquad \leq E\bigl[|\xi|^2e^{\int_0^Td(s)\,ds}\bigr]+C_{\lambda}\int_t^T
E\bigl[e^{\int_0^sd(u)\,du}Y_n^2(s)\bigr]\,ds\\
&&\qquad \quad {}+E\biggl[\int_t^T e^{\int_0^sd(u)\,du}|f(s,0,0)|^2\,ds\biggr].\nonumber
\end{eqnarray}
Gronwall's inequality yields
%
\begin{eqnarray} \label{2.8}
&&\sup_n\sup_{0\leq t\leq T}E\bigl[|Y_n(t)|^2e^{\int_0^td(s)\,ds}\bigr]\nonumber\\[-8pt]\\[-8pt]
&&\qquad \leq
C\biggl\{E\bigl[|\xi|^2e^{\int_0^Td(s)\,ds}\bigr]+E\biggl[\int_0^T
e^{\int_0^sd(u)\,du}|f(s,0,0)|^2\,ds\biggr]\biggr\}\nonumber
\end{eqnarray}
and also
%
\begin{equation}\label{2.9}
\sup_nE\biggl[\int_0^Te^{\int_0^sd(u)\,du}\langle a(X(s))Z_n(s), Z_n(s)\rangle\,ds
\biggr]<\infty.
\end{equation}
Moreover, (\ref{2.6})--(\ref{2.9}) further imply that there exists
some constant $C$ such~that
%
\begin{eqnarray}\label{2.10}
&& E\Bigl[\sup_{0\leq t\leq T}Y_n^2(t)e^{\int_0^td(s)\,ds}\Bigr]\nonumber\\
&&\qquad \leq C+CE\biggl[\sup_{0\leq t\leq
T}\int_0^te^{\int_0^sd(u)\,du}Y_n(s)\langle Z_n(s),dM(s)\rangle\biggr]\nonumber\\
&&\qquad \leq C+CE \biggl[ \biggl(
\int_0^Te^{2\int_0^sd(u)\,du}Y_n^2(s)\langle a(X(s))Z_n(s), Z_n(s)\rangle\,ds
\biggr)^{\gfrac{1}{2}} \biggr]\nonumber\\
&&\qquad \leq C+CE \biggl[\sup_{0\leq s\leq
T}\bigl(e^{\vfrac{1}{2}\int_0^sd(u)\,du}|Y_n(s)|\bigr)\nonumber\\
&&\qquad \quad \hphantom{C+CE \biggl[}
{}\times \biggl(
\int_0^Te^{\int_0^sd(u)\,du}\langle a(X(s))Z_n(s), Z_n(s)\rangle\,ds
\biggr)^{\gfrac{1}{2}} \biggr]\\
&&\qquad \leq C+\frac{1}{2}E \Bigl[\sup_{0\leq s\leq
T}\bigl(e^{\int_0^sd(u)\,du}Y_n^2(s)\bigr) \Bigr]\nonumber\\
&&\qquad \quad{} +C_1E \biggl[ \int_0^Te^{\int_0^sd(u)\,du}\langle a(X(s))Z_n(s),
Z_n(s)\rangle\,ds \biggr].\nonumber
\end{eqnarray}
In view of (\ref{2.9}), this yields
%
\begin{equation}\label{2.11}
\sup_nE\Bigl[\sup_{0\leq t\leq T}Y_n^2(t)e^{\int_0^td(s)\,ds}\Bigr]<\infty.
\end{equation}
By (\ref{2.9}) and (\ref{2.11}), we can extract a subsequence $n_k$
such that $Y_{n_k}(t)e^{\vfrac{1}{2}\int_0^td(s)\,ds}$ converges to
some $\hat{Y}(t)$ in $L^2(\Omega, L^{\infty}[0,T])$ equipped with
the weak star topology and $Z_{n_k}(t)e^{\vfrac{1}{2}\int_0^td(s)\,ds}$
converges weakly to some $\hat{Z}(t)$ in $L^2(\Omega_T; R)$, where
$\Omega_T=[0,T]\times\Omega$. Observe that
%
\begin{eqnarray}\label{2.12}
&&Y_{n_k}(t)e^{\vfrac{1}{2}\int_0^td(s)\,ds}\nonumber\\
&&\qquad =e^{\vfrac{1}{2}\int_0^Td(s)\,ds}\xi+\int_t^T
e^{\vfrac{1}{2}\int_0^sd(u)\,du}f_{n_k}(s,Y_{n_k}(s),
Z_{n_k}(s))\,ds\nonumber\\[-8pt]\\[-8pt]
&&\qquad \quad {}-\frac{1}{2}\int_t^Te^{\vfrac{1}{2}\int_0^sd(u)\,du}Y_{n_k}(s)d(s)\,ds\nonumber \\
&&\qquad \quad {}-
\int_t^Te^{\vfrac{1}{2}\int_0^sd(u)\,du}\langle Z_{n_k}(s),dM(s)\rangle.\nonumber
\end{eqnarray}
Letting $k\rightarrow\infty$ in (\ref{2.12}), using the
monotonicity of $f$, following the same arguments as that in the proof
of Proposition 2.3 in
Darling and Pardoux in \cite{DP}, we can show that the limit $(\hat{Y},
\hat{Z})$ satisfies
%
\begin{eqnarray}\label{2.13}
\hat{Y}(t)&=&e^{\vfrac{1}{2}\int_0^Td(s)\,ds}\xi\nonumber\\
&&{}+\int_t^T
e^{\vfrac{1}{2}\int_0^sd(u)\,du}f\bigl(s,e^{-\vfrac{1}{2}\int_0^sd(u)\,du}\hat
{Y}(s), e^{-\vfrac{1}{2}\int_0^sd(u)\,du}\hat{Z}(s)\bigr)\,ds\hspace*{-10pt}\\
&&{}-\frac{1}{2}\int_t^T\hat{Y}(s)d(s)\,ds- \int_t^T\langle \hat{Z}(s),dM(s)\rangle.\nonumber
\end{eqnarray}
Set
\[
Y(t)=e^{-\vfrac{1}{2}\int_0^td(u)\,du}\hat{Y}(t),\qquad
Z(t)=e^{-\vfrac{1}{2}\int_0^td(u)\,du}\hat{Z}(t).
\]
An application of
It\^o's formula yields that
\[
Y(t)=\xi+\int_t^Tf(s,Y(s), Z(s))\,ds- \int_t^T\langle Z(s),dM(s)\rangle,
\]
namely, $(Y, Z)$ is a solution to the backward equation (\ref{2.1}).
The proof is complete.
\end{pf}

\subsection{BSDEs with random terminal times}\label{sec3.2}
Let $f(t,y,z)$ satisfy (A.1)--(A.3) in Section~\ref{sec3.1}. In this
subsection, set $d(s)=-2d_1(s)+\delta d_2^2$. The following result
provides existence and uniqueness for BSDEs with random terminal
time. Let $\tau$ be a stopping time. Suppose $\xi$ is $\mathcal{
F}_{\tau}$-measurable.
\begin{theorem}\label{thm3.2}
Assume $E[e^{\int_0^{\tau}d(s)\,ds}|\xi|^2]<\infty$,  $E[\int_0^\tau K(s)\,ds]<\infty$ and
%
\begin{equation}\label{2.14}
E \biggl[\int_0^{\tau}e^{\int_0^{s}d(u)\,du}|f(s,0,0)|^2\,ds
\biggr]<\infty,
\end{equation}
for some $\delta
>\frac{1}{\lambda}$, where $\lambda$ is the constant appeared in (2.1).
Then, there exists a
unique solution $(Y, Z)$ to the BSDE
%
\begin{equation}\label{2.15}
Y(t)=\xi+\int_{\tau\wedge t}^{\tau}f(s, Y(s),
Z(s))\,ds-\int_{\tau\wedge t}^{\tau}\langle Z(s), dM(s)\rangle.
\end{equation}
Furthermore, the solution $(Y, Z)$ satisfies
%
\begin{equation}\label{2.16}
 E \biggl[\int_0^{\tau}e^{\int_0^{s}d(u)\,du}Y^2(s)\,ds \biggr]<\infty,\qquad
E \biggl[\int_0^{\tau}e^{\int_0^{s}d(u)\,du}|Z(s)|^2\,ds
\biggr]<\infty,\hspace*{-35pt}
\end{equation}
and
%
\begin{equation}\label{2.17}
E \Bigl[ \sup_{0\leq s\leq\tau}\bigl\{e^{\int_0^{s}d(u)\,du}Y^2(s)\bigr\}
\Bigr]<\infty.
\end{equation}
\end{theorem}

\begin{pf}
After the preparation of Theorem
\ref{thm3.1}, the proof of this theorem is similar to that of Theorem 3.4 in
\cite{DP}, where $d_1(s)$, $d_2$ were both assumed to be constants. We
only give a sketch highlighting the differences. For every $n\geq
1$, from Theorem~\ref{thm3.1} we know that the following BSDE has a unique
solution $(Y_n, Z_n)$ on $0\leq t\leq n$:
%
\begin{equation}\label{2.18}
  Y_n(t)=E[\xi|\mathcal{ F}_n]+\int_{\tau\wedge t}^{\tau\wedge n}f(s,
Y_n(s), Z_n(s))\,ds-\int_{\tau\wedge t}^{\tau\wedge n}\langle Z_n(s), dM(s)\rangle.\hspace*{-35pt}
\end{equation}
Extend the definition of $(Y_n, Z_n)$ to all $t\geq0$ by setting
\[
Y_n(t)=E[\xi|\mathcal{ F}_n],\qquad  Z_n(t)=0\qquad  \mbox{for
} t\geq n.
\]
Then the extended $(Y_n, Z_n)$ satisfies a bsde
similar to (\ref{2.18}) with $f$ replaced by $\chi_{\{s\leq n\wedge
\tau\}}f(s, y,z)$. Let $n\geq m$. By It\^o's formula, we have
%
\begin{eqnarray}\label{2.19}
&&|Y_n(t\wedge\tau)-Y_m(t\wedge\tau)|^2e^{\int_0^{t\wedge
\tau}d(s)\,ds}\nonumber\hspace*{-35pt}\\
&&\quad {}+\int_{t\wedge\tau}^{n\wedge
\tau}e^{\int_0^sd(u)\,du}\bigl\langle a(X(s))\bigl(Z_n(s)-Z_m(s)\chi_{\{s\leq m\wedge
\tau\}}\bigr),\nonumber \hspace*{-35pt}\\
&&\hspace*{144pt}
 Z_n(s)-Z_m(s)\chi_{\{s\leq m\wedge
\tau\}}\bigr\rangle\,ds\nonumber \hspace*{-35pt}\\
&&\qquad =e^{\int_0^{n\wedge\tau}d(s)\,ds} (E[\xi|\mathcal{ F}_n]-E[\xi
|\mathcal{ F}_m] )^2\hspace*{-35pt}\\
&&\qquad \quad {}-\int_{t\wedge\tau}^{n\wedge\tau} e^{\int_0^{s\wedge
\tau}d(u)\,du}|Y_n(s\wedge\tau)-Y_m(s\wedge\tau)|^2d(s)\,ds\nonumber\hspace*{-35pt}\\
&&\qquad \quad {} +2\int_{t\wedge
\tau}^{n\wedge\tau}e^{\int_0^{s\wedge\tau}d(u)\,du}\bigl(Y_n(s\wedge
\tau)-Y_m(s\wedge\tau)\bigr)\nonumber\hspace*{-35pt}\\
&&\qquad \quad\hphantom{{} +2\int_{t\wedge\tau}^{n\wedge\tau}}
{} \times \bigl( f\bigl(s,Y_n(s\wedge\tau), Z_n(s\wedge
\tau)\bigr)\nonumber\hspace*{-35pt} \\
&&\qquad \quad\hphantom{{} +2\int_{t\wedge\tau}^{n\wedge\tau}\times \bigl( }
{}-f\bigl(s,Y_m(s\wedge\tau), Z_m(s\wedge\tau)\bigr)
\bigr)\,ds\nonumber\hspace*{-35pt}\\
&&\qquad \quad {}+2\int_{m\wedge
\tau}^{n\wedge\tau}e^{\int_0^{s\wedge\tau}d(u)\,du}\bigl(Y_n(s\wedge
\tau)-Y_m(s\wedge\tau)\bigr)\nonumber \hspace*{-35pt}\\
&&\quad \qquad\hphantom{{}+2\int_{m\wedge\tau}^{n\wedge\tau}}
{}\times f\bigl(s,Y_m(s\wedge\tau), Z_m(s\wedge
\tau)\bigr)\,ds\nonumber\hspace*{-35pt}\\
&&\qquad \quad {}-2\int_{t\wedge
\tau}^{n\wedge\tau}e^{\int_0^{s\wedge\tau}d(u)\,du}\bigl(Y_n(s\wedge
\tau)-Y_m(s\wedge\tau)\bigr)\langle Z_n(s\wedge\tau), dM(s)\rangle\nonumber\hspace*{-35pt}\\
&&\qquad \quad {}+2\int_{t\wedge
\tau}^{m\wedge\tau}e^{\int_0^{s\wedge\tau}d(u)\,du}\bigl(Y_n(s\wedge
\tau)-Y_m(s\wedge\tau)\bigr)\langle Z_m(s\wedge\tau), dM(s)\rangle.\nonumber\hspace*{-35pt}
\end{eqnarray}
Choose $\delta_1$, $\delta_2$ such that
$\frac{1}{\lambda}<\delta_1<\delta$ and
$0<\delta_2<\delta-\delta_1$. In view of the (A.1) and (A.2), we
have
%
\begin{eqnarray}\label{2.20}
&& 2\int_{t\wedge
\tau}^{n\wedge\tau}e^{\int_0^{s\wedge\tau}d(u)\,du}\bigl(Y_n(s\wedge
\tau)-Y_m(s\wedge\tau)\bigr)\nonumber\\
&&\hphantom{2\int_{t\wedge\tau}^{n\wedge\tau}}
{} \times \bigl( f\bigl(s,Y_n(s\wedge\tau), Z_n(s\wedge
\tau)\bigr)-f\bigl(s,Y_m(s\wedge\tau), Z_m(s\wedge\tau)\bigr)
\bigr)\,ds\nonumber\\
&&\qquad \leq-2\int_{t\wedge
\tau}^{n\wedge\tau}e^{\int_0^{s\wedge\tau}d(u)\,du}\bigl(Y_n(s\wedge
\tau)-Y_m(s\wedge\tau)\bigr)^2d_1(s)\,ds\nonumber\\[-8pt]\\[-8pt]
&&\qquad \quad {}+\delta_1d_2^2\int_{t\wedge
\tau}^{n\wedge\tau}e^{\int_0^{s\wedge\tau}d(u)\,du}\bigl(Y_n(s\wedge
\tau)-Y_m(s\wedge\tau)\bigr)^2\,ds\nonumber\\
&&\qquad \quad {}+\frac{1}{\lambda\delta_1}\int_{t\wedge\tau}^{n\wedge
\tau}e^{\int_0^sd(u)\,du}\bigl\langle a(X(s))\bigl(Z_n(s)-Z_m(s)\chi_{\{s\leq m\wedge
\tau\}}\bigr),\nonumber\\
&&\hspace*{186pt} Z_n(s)-Z_m(s)\chi_{\{s\leq m\wedge\tau\}}\bigr\rangle\,ds.\nonumber
\end{eqnarray}
On the other hand, by (A.3), it follows that
%
\begin{eqnarray}\label{2.21}
&& 2\int_{m\wedge
\tau}^{n\wedge\tau}e^{\int_0^{s\wedge\tau}d(u)\,du}\bigl(Y_n(s\wedge
\tau)-Y_m(s\wedge\tau)\bigr)f\bigl(s,Y_m(s\wedge\tau), Z_m(s\wedge
\tau)\bigr)\,ds\nonumber\hspace*{-35pt}\\
&&\qquad \leq \delta_2d_2^2\int_{t\wedge\tau}^{n\wedge\tau}
e^{\int_0^{s\wedge\tau}d(u)\,du}\bigl(Y_n(s\wedge\tau)-Y_m(s\wedge
\tau)\bigr)^2\,ds\hspace*{-35pt}\\
&&\qquad \quad {}+\frac{1}{\delta_2 d_2^2}\int_{m\wedge
\tau}^{n\wedge\tau}e^{\int_0^{s\wedge\tau}d(u)\,du}\bigl(
|f(s,0,0)|+K+K|E[\xi|\mathcal{
F}_m]|\bigr)^2\,ds.\nonumber\hspace*{-35pt}
\end{eqnarray}

\noindent
Take expectation and utilize (\ref{2.19})--(\ref{2.21}) to obtain
%
\begin{eqnarray}\label{2.22}
&&E \bigl[|Y_n(t\wedge\tau)-Y_m(t\wedge\tau)|^2e^{\int_0^{t\wedge
\tau}d(s)\,ds} \bigr]\nonumber\hspace*{-35pt}\\
&&\quad {}+\biggl(1-\frac{1}{\lambda\delta_1}\biggr)E \biggl[\int_{t\wedge\tau}^{n\wedge
\tau}e^{\int_0^sd(u)\,du}\bigl\langle a(X(s))\bigl(Z_n(s)-Z_m(s)\chi_{\{s\leq m\wedge
\tau\}}\bigr) , \nonumber\hspace*{-35pt}\\
&&\hspace*{207pt} Z_n(s)-Z_m(s)\chi_{\{s\leq m\wedge
\tau\}}\bigr\rangle\,ds \biggr] \nonumber\hspace*{-35pt}\\[-8pt]\\[-8pt]
&&\quad {}+(\delta-\delta_1-\delta_2)d_2^2E
\biggl[\int_{t\wedge\tau}^{n\wedge\tau}e^{\int_0^sd(u)\,du}\bigl(Y_n(s\wedge
\tau)-Y_m(s\wedge
\tau)\bigr)^2\,ds \biggr]\nonumber\hspace*{-35pt}\\
&&\qquad \leq E \bigl[e^{\int_0^{n\wedge\tau}d(s)\,ds} (E[\xi|\mathcal{
F}_n]-E[\xi|\mathcal{ F}_m] )^2 \bigr]\nonumber\hspace*{-35pt}\\
&&\qquad \quad {}+\frac{1}{\delta_2 d_2^2}E \biggl[ \int_{m\wedge
\tau}^{n\wedge\tau}e^{\int_0^{s\wedge\tau}d(u)\,du}\bigl(
|f(s,0,0)|+K+K|E[\xi|\mathcal{ F}_m]|\bigr)^2\,ds \biggr].\nonumber\hspace*{-35pt}
\end{eqnarray}
Since the right-hand side tends to zero as $n,m\rightarrow\infty$, we
deduce that
\[
\bigl\{\bigl(e^{\vfrac{1}{2}\int_0^{t\wedge
\tau}d(s)\,ds}Y_n(t),e^{\vfrac{1}{2}\int_0^{t\wedge
\tau}d(s)\,ds}Z_n(t)\bigr)\bigr\}
\]
converges to some $(\hat{Y}, \hat{Z})$ in
$M^2(0, \tau; R\times R^d)$. Furthermore, for every $t\geq0$,
$e^{\vfrac{1}{2}\int_0^{t\wedge\tau}d(s)\,ds}Y_n(t)$ converges in
$L^2$. We may as well assume
%
\begin{equation}\label{2.23}
\hat{Y}(t)=\lim_{n\rightarrow\infty}e^{\vfrac{1}{2}\int_0^{t\wedge
\tau}d(s)\,ds}Y_n(t)
\end{equation}
for all $t$. Observe that for any $n\geq t\geq0$,
%
\begin{eqnarray}\label{2.24}
&&e^{\vfrac{1}{2}\int_0^{t\wedge
\tau}d(s)\,ds}Y_n(t)\nonumber\hspace*{-35pt}\\
&&\qquad =e^{\vfrac{1}{2}\int_0^{n\wedge\tau}d(s)\,ds}E[\xi|\mathcal{ F}_n]
+\int_{\tau\wedge t}^{n\wedge\tau}e^{\vfrac{1}{2}\int_0^{s\wedge
\tau}d(u)\,du}f(s, Y_n(s),
Z_n(s))\,ds\nonumber\hspace*{-35pt}\\[-8pt]\\[-8pt]
&&\qquad \quad {}-\frac{1}{2}\int_{\tau\wedge t}^{n\wedge
\tau}e^{\vfrac{1}{2}\int_0^{s\wedge\tau}d(u)\,du}
Y_n(s)d(s)\,ds\nonumber\hspace*{-35pt}\\
&&\qquad \quad {}-\int_{\tau\wedge t}^{n\wedge\tau}e^{\vfrac{1}{2}\int_0^{s\wedge
\tau}d(u)\,du}\langle Z_n(s), dM(s)\rangle.\nonumber\hspace*{-35pt}
\end{eqnarray}
Letting $n\rightarrow\infty$ yields that
%
\begin{eqnarray}\label{2.25}
\hat{Y}(t)&=&e^{\vfrac{1}{2}\int_0^{\tau}d(s)\,ds}\xi+\int_{\tau\wedge
t}^{\tau}e^{\vfrac{1}{2}\int_0^{s\wedge\tau}d(u)\,du}f\bigl(s,
e^{-\vfrac{1}{2}\int_0^{s\wedge\tau}d(u)\,du}\hat{Y}(s),\nonumber\hspace*{-35pt}\\
&&\hspace*{195pt}
e^{-\vfrac{1}{2}\int_0^{s\wedge
\tau}d(u)\,du}\hat{Z}(s)\bigr)\,ds\hspace*{-35pt}\\
&&{}-\frac{1}{2}\int_{\tau\wedge t}^{
\tau}\hat{Y}(s)d(s)\,ds-\int_{\tau\wedge t}^{\tau}\langle\hat{Z}(s), dM(s)\rangle.\nonumber\hspace*{-35pt}
\end{eqnarray}
Put
\[
Y(t)=e^{-\vfrac{1}{2}\int_0^{t\wedge
\tau}d(s)\,ds}\hat{Y}(t),\qquad  Z(t)=e^{-\vfrac{1}{2}\int_0^{t\wedge
\tau}d(s)\,ds}\hat{Z}(t).
\]
An application of It\^o's formula and
(\ref{2.25}) yield that

\begin{equation}\label{2.26}
Y(t)=\xi+\int_{\tau\wedge t}^{\tau}f(s, Y(s),
Z(s))\,ds-\int_{\tau\wedge t}^{\tau}\langle Z(s), dM(s)\rangle.
\end{equation}
Hence, $(Y, Z)$ is a solution to the bsde (\ref{2.15}) proving the
existence. To obtain the estimates (\ref{2.16}) and (\ref{2.17}), we
proceed to get an uniform estimate for $Y_n(s)$ and then pass to the
limit. Let $\delta_1, \delta_2$ be chosen as before. Similar to the
proof of (\ref{2.8}), by It\^o's formula, we have
%
\begin{eqnarray}\label{2.27}
&&|Y_n(t\wedge\tau)|^2e^{\int_0^{t\wedge\tau}d(s)\,ds}+\int_{t\wedge
\tau}^{n\wedge\tau} e^{\int_0^sd(u)\,du}\langle a(X(s))Z_n(s),
Z_n(s)\rangle\,ds \nonumber\\
&&\qquad \leq|E[\xi|\mathcal{ F}_n]|^2e^{\int_0^{n\wedge
\tau}d(s)\,ds}-\int_{t\wedge\tau}^{n\wedge
\tau} e^{\int_0^sd(u)\,du}|Y_n(s)|^2d(s)\,ds\nonumber\\
&&\qquad \quad {} -2\int_{t\wedge\tau}^{n\wedge\tau}
e^{\int_0^sd(u)\,du}d_1(s)Y_n^2(s)\,ds\nonumber\\
&&\qquad \quad {}+2\int_{t\wedge\tau}^{n\wedge
\tau} e^{\int_0^sd(u)\,du}d_2|Y_n(s)||Z_n(s)|\,ds\nonumber\\
&&\qquad \quad {}+2\int_{t\wedge\tau}^{n\wedge\tau}
e^{\int_0^sd(u)\,du}|Y_n(s)||f(s,0,0)|\,ds\nonumber\\
&&\qquad \quad {}-2\int_{t\wedge\tau}^{n\wedge
\tau} e^{\int_0^sd(u)\,du}Y_n(s)\langle
Z_n(s),dM(s)\rangle\\
&&\qquad \leq |E[\xi|\mathcal{ F}_n]|^2e^{\int_0^{n\wedge
\tau}d(s)\,ds}-\int_{t\wedge\tau}^{n\wedge\tau} e^{\int
_0^sd(u)\,du}\delta d_2^2Y_n^2(s)\,ds\nonumber\\
&&\qquad \quad {}+\int_{t\wedge\tau}^{n\wedge\tau}
e^{\int_0^sd(u)\,du}\delta_1 d_2^2Y_n^2(s)\,ds\nonumber\\
&&\qquad \quad {}+\frac{1}{\delta_1\lambda}\int_{t\wedge\tau}^{n\wedge\tau}
e^{\int_0^sd(u)\,du}\langle a(X(s))Z_n(s), Z_n(s)\rangle\, ds
\nonumber\\
&&\qquad \quad {}+\int_{t\wedge\tau}^{n\wedge\tau}\delta_2d_2^2
e^{\int_0^sd(u)\,du}Y_n^2(s)\,ds+\frac{1}{\delta_2 d_2^2}\int_{t\wedge
\tau}^{n\wedge\tau}
e^{\int_0^sd(u)\,du}|f(s,0,0)|^2\,ds\hspace*{-8pt}\nonumber\\
&&\qquad \quad {}-2\int_{t\wedge\tau}^{n\wedge\tau}
e^{\int_0^sd(u)\,du}Y_n(s)\langle Z_n(s),dM(s)\rangle.\nonumber
\end{eqnarray}
Recalling the choices of $d(s)$, $\delta_1$ and $\delta_2$, using
Burkholder's inequality, we obtain from (\ref{2.27}) that
%
\begin{eqnarray}\label{2.28}
\qquad &&E \Bigl[ \sup_{0\leq t\leq n}|Y_n(t\wedge
\tau)|^2e^{\int_0^{t\wedge\tau}d(s)\,ds} \Bigr]\nonumber\\
\qquad &&\qquad \leq E\bigl[|\xi|^2e^{\int_0^{n\wedge\tau}d(s)\,ds}\bigr]+E \biggl[\int_{0}^{
\tau}e^{\int_0^sd(u)\,du}\frac{1}{\delta_2 d_2^2}|f(s,0,0)|^2\,ds \biggr]\\
\qquad &&\qquad \quad {}+2CE \biggl[ \biggl(\int_{0}^{n\wedge\tau}
e^{2\int_0^sd(u)\,du}Y_n^2(s)\langle a(X(s))Z_n(s), Z_n(s)\rangle\,ds \biggr)^{\gfrac{1}{2}}
\biggr]\nonumber\\
\qquad &&\qquad \leq E\bigl[|\xi|^2e^{\int_0^{n\wedge\tau}d(s)\,ds}\bigr]+E \biggl[\int_{0}^{
\tau}
e^{\int_0^sd(u)\,du}\frac{1}{\delta_2 d_2^2}|f(s,0,0)|^2\,ds \biggr]\nonumber\\
\qquad &&\qquad \quad {}+\frac{1}{2}E \Bigl[ \sup_{0\leq t\leq n}|Y_n(t\wedge
\tau)|^2e^{\int_0^{t\wedge\tau}d(s)\,ds} \Bigr]\nonumber\\
\qquad &&\qquad \quad {}+C_1E \biggl[ \int_{0}^{n\wedge\tau}
e^{\int_0^sd(u)\,du}\langle a(X(s))Z_n(s), Z_n(s)\rangle\,ds \biggr].\nonumber
\end{eqnarray}
In view of (\ref{2.27}), as the proof of (\ref{2.9}), we can show
that
%
\begin{equation}\label{2.29}
\sup_nE \biggl[ \int_{0}^{n\wedge\tau}
e^{\int_0^sd(u)\,du}\langle a(X(s))Z_n(s), Z_n(s)\rangle\,ds \biggr]<\infty.
\end{equation}
(\ref{2.29}) and (\ref{2.28}) together with our assumptions on $f$
and $\xi$ imply
%
\begin{equation}\label{2.30}
\sup_n E \Bigl[ \sup_{0\leq t\leq n}|Y_n(t\wedge
\tau)|^2e^{\int_0^{t\wedge\tau}d(s)\,ds} \Bigr]<\infty.
\end{equation}
Applying Fatou lemma, (\ref{2.17}) follows.
\end{pf}

\subsection{A particular case}\label{sec3.3}
Let $f(x,y,z)\dvtx R^d\times R\times R^d\rightarrow R$ be a Borel
measurable function. Assume that $f$ is continuous in $y$ and
satisfies:
\begin{itemize}[(B.1)]
\item[(B.1)] $(y_1-y_2)(f(x,y_1,z)-f(x,y_2,z))\leq-c_1(x)|y_1-y_2|^2$,
where $c_1(x)$ is a measurable function
on $R^d$.
\item[(B.2)] $|f(x,y,z_1)-f(x,y,z_2)|\leq c_2 |z_1-z_2|$.
\item[(B.3)] $|f(x,y,z)|\leq|f(x,0,z)|+c_3(x)(1+|y|).$
\end{itemize}
Let $D$ be a bounded regular domain. Define
%
\begin{equation}\label{2.31}
\tau_D^x=\inf\{t\geq0\dvtx X_x(t)\notin D\}.
\end{equation}
Given $g\in C_b(R^d)$. Consider for each $x\in D$ the following
BSDE:
%
\begin{eqnarray}\label{2.32}
Y_x(t)&=& g(X_x(\tau_D^x))+\int_{t\wedge\tau_D^x}^{\tau_D^x}
f(X_x(s),Y_x(s), Z_x(s))\,ds\nonumber\\[-8pt]\\[-8pt]
&&{}-\int_{t\wedge\tau_D^x}^{\tau_D^x}\langle Z_x(s),dM_x(s)\rangle,\nonumber
\end{eqnarray}
where $M_x(s)$ is the martingale part of $X_x(s)$. As a consequence
of Theorem~\ref{thm3.2}, we have
the following theorem.
\begin{theorem}\label{thm3.3}
Suppose $c_3\in L^P(D)$ for $p>\frac{d}{2}$,
\[
E_x\biggl[\exp\biggl(\int_0^{\tau_D^x} \bigl(-2c_1(X(s))+\delta c_2^2 \bigr)\,ds\biggr)\biggr]<\infty,
\]
for some $\delta>\frac{1}{\lambda}$ and
\[
E_x\biggl[\int_0^{\tau_D^x}|f(X(s),0,0)|^2\,ds\biggr]<\infty.
\]
The BSDE (\ref{2.32}) admits a unique solution $(Y_x(t), Z_x(t))$.
Furthermore,
%
\begin{equation}\label{2.33}
\sup_{x\in\bar{D}}|Y_x(0)|<\infty.
\end{equation}
\end{theorem}
%
\section{Semilinear PDEs}\label{sec4}

As in previous sections, $(X(t), P_x)$ will denote the diffusion
process defined in (\ref{1.3}).
\subsection{Linear case}\label{sec4.1}
Consider the
second order differential operator
%
\begin{equation}\label{3.1}
L_2={1\over2}\sum_{i,j=1}^{d}{\partial\over{\partial
x_i}}\biggl(a_{ij}(x)\,{\partial\over{\partial
x_j}}\biggr)+\sum_{i=1}^{d}b_i(x)\,{\partial\over{\partial
x_i}}+q(x).\vadjust{\goodbreak}
\end{equation}
Let $D$ be a bounded domain with regular boundary (w.r.t. the
Laplace operator
$\Delta$) and $F(x)$ a measurable function satisfying
%
\begin{equation}\label{3.2}
|F(x)|\leq C+C|q(x)|.
\end{equation}
Take $\varphi\in C(\partial D)$ and consider the Dirichlet boundary
value problem
%
\begin{equation}\label{3.3}
\cases{
L_2 u= F, &\quad in $D$, \cr
u= \varphi,&\quad on $\partial D$.
}
\end{equation}

\begin{theorem}\label{thm4.1}
Assume (\ref{3.2}) and that there exists $x_0\in D$ such that
\[
E_{x_0}\biggl[\exp\biggl(\int_0^{\tau_D^{x_0}}q(X(s))\,ds\biggr)\biggr]<\infty.
\]
Then there is a
unique, continuous weak solution $u$ to the Dirichlet boundary value
problem (\ref{3.3}) which is given by
%
\begin{equation}\label{3.4}
u(x)=E_x\biggl[\varphi(X(\tau_D^x))+\int_0^{\tau_D^x}
e^{\int_0^tq(X(s))\,ds}F(X(t)\,dt\biggr].
\end{equation}
\end{theorem}

\begin{pf}
Write
\[
u_1(x)=E_x[\varphi(X(\tau_D^x))],
\]
and
\[
u_2(x)=E_x \biggl[\int_0^{\tau_D^x}e^{\int_0^tq(X(s))\,ds}F(X(t))\,dt
\biggr].
\]
We know from Theorem 4.3 in \cite{CZh} that $u_1$ is the unique,
continuous weak solution to the problem
%
\begin{equation}\label{3.5}
\cases{
 L_2 u= 0, &\quad in $D$, \cr
u= \varphi,&\quad on $\partial D$.
}
\end{equation}
So it is sufficient to show that $u_2$ is the unique,
continuous weak solution to the following problem:
%
\begin{equation}\label{3.5}
\cases{
L_2 u= F, &\quad in $ D$, \cr
u= 0, &\quad on $\partial D$.
}
\end{equation}
By Lemma 5.7 in \cite{CZh} and Proposition 3.16 in \cite{ChungZ},
we know that $u_2$ belong to $C_0(D)$. Let $G_{\beta}, \beta\geq0$
denote the
resolvent operators of the generator $L_2$ on $D$ with Dirichlet
boundary condition, that is,
\[
G_{\beta}f(x)=E_x \biggl[\int_0^{\tau_D^x}e^{-\beta t}e^{\int
_0^tq(X(s))\,ds}f(X(t))\,dt \biggr].
\]
By the Markov property, it is easy to see that
\[
\beta\bigl(u_2(x)-\beta G_{\beta}u_2(x)\bigr)=\beta G_{\beta}F(x).\vadjust{\goodbreak}
\]
Since $G_{\beta}$ is strong continuous, it follows that
\[
\lim_{\beta\rightarrow\infty}\beta(u_2-\beta G_{\beta}u_2)=F
\]
in $L^2(D)$. This shows that $u_2\in D(L_2)\subset W^{1,2}(D)$ and
$ L_2u_2=F$. The proof is complete.
\end{pf}

\subsection{Semilinear case}\label{sec4.2}
Let $g(x,y,z)\dvtx R^d\times R\times R^d\rightarrow R$ be a Borel
measurable function that satisfies:
\begin{itemize}[(C.1)]
\item[(C.1)] $(y_1-y_2)(g(x,y_1,z)-g(x,y_2,z))\leq-k_1(x)|y_1-y_2|^2,$
\item[(C.2)] $|g(x,y,z_1)-g(x,y,z_2)|\leq k_2 |z_1-z_2|$,
\item[(C.3)] $|g(x,y,z)|\leq C+C|q(x)|,$
\end{itemize}
where $k_1(x)$ is a measurable function and $k_2, C$ are constants.
Consider the semilinear Dirichlet boundary value problem
%
\begin{equation}\label{3.6}
\cases{
L_2 u= -g(x, u(x), \nabla u(x)), &\quad in $D$, \cr
u= \varphi,&\quad on $\partial D$,
}
\end{equation}
where $\varphi\in C(\partial D)$.
\begin{theorem}\label{thm4.2}
Assume
\[
E_x\biggl[\exp\biggl(\int_0^{\tau_D^x} \bigl(q(X(s))-2k_1(X(s))+\delta k_2^2 \bigr)\,ds\biggr)\biggr]<\infty,
\]
for some $\delta>\frac{1}{\lambda}$. 

The Dirichlet boundary value problem (\ref{3.6}) has a unique
continuous weak solution.
\end{theorem}

\begin{pf}
Set $f(x,y,z)=q(x)y+g(x,y,z)$. According to Theorem
\ref{thm3.3},
for every $x\in D$ the following BSDE:
%
\begin{eqnarray}\label{3.7}
Y_x(t)&=& \phi(X_x(\tau_D^x))+\int_{t\wedge\tau_D^x}^{\tau_D^x}
f(X_x(s),Y_x(s), Z_x(s))\,ds\nonumber\\[-8pt]\\[-8pt]
&&{}-\int_{t\wedge\tau_D^x}^{\tau_D^x}\langle Z_x(s),dM_x(s)\rangle,\nonumber
\end{eqnarray}
admits a unique solution $(Y_x(t), Z_x(t)), t\geq0$. Put
$u_0(x)=Y_x(0)$ and $v_0(x)=Z_x(0)$. By the strong Markov property
of $X$ and the uniqueness of the BSDE (\ref{3.7}), it is easy to see
that
%
\begin{equation}
Y_x(t)=u_0(X_x(t)),\qquad  Z_x(t)=v_0(X_x(t)),\qquad  0\leq t\leq
\tau_D^x.
\end{equation}
Now consider the following problem:
%
\begin{equation}\label{3.8}
\cases{
L_1 u= -f(x, u_0(x), v_0(x)), &\quad in $D$, \cr
u= \varphi,&\quad on $\partial D$,
}
\end{equation}
where $L_1$ is defined as in Section~\ref{sec2}.
By Theorem~\ref{thm4.1}, problem (\ref{3.8}) has a unique continuous weak solution
$u(x)$. As $u\in W^{1,2}(D)$, it follows from the decomposition of
the Dirichlet process $u(X(t\wedge\tau_D^x))$ (see \cite{FOT,LZ}) that
%
\begin{eqnarray}\label{3.9}
u\bigl(X(t\wedge\tau_D^x)\bigr)&=& \varphi(X_x(\tau_D^x))+\int_{t\wedge
\tau_D^x}^{\tau_D^x}
f\bigl(X_x(s),u_0\bigl(X(s\wedge\tau_D^x)\bigr), v_0\bigl(X(s\wedge\tau_D^x)\bigr)\bigr)\,ds\nonumber\\
&&{}-\int_{t\wedge\tau_D^x}^{\tau_D^x}\bigl\langle\nabla u\bigl(X(s\wedge
\tau_D^x)\bigr),dM_x(s)\bigr\rangle\nonumber\\[-8pt]\\[-8pt]
&=& \varphi(X_x(\tau_D^x))+\int_{t\wedge\tau_D^x}^{\tau_D^x}
f(X_x(s),Y_x(s), Z_x(s)))\,ds\nonumber\\
&&{}-\int_{t\wedge\tau_D^x}^{\tau_D^x}\bigl\langle\nabla u\bigl(X(s\wedge
\tau_D^x)\bigr),dM_x(s)\bigr\rangle.\nonumber
\end{eqnarray}
Take conditional expectation both in (\ref{3.9}) and (\ref{3.7}) to
discover
\begin{eqnarray*}
Y_x(t\wedge\tau_D^x)&=&u\bigl(X(t\wedge\tau_D^x)\bigr)\\
&=&E \biggl[\varphi(X_x(\tau
_D^x))+\int_{t\wedge\tau_D^x}^{\tau_D^x}
f(X_x(s),Y_x(s), Z_x(s))\,ds \Big|\mathcal{ F}_{t\wedge\tau_D^x} \biggr].
\end{eqnarray*}
In
particular, let $t=0$ to obtain $u(x)=u_0(x)$. On the other hand,
comparing (\ref{3.7}) with (\ref{3.9}) and by the uniqueness of
decomposition of semimartingales, we deduce that
\[
\int_{t\wedge\tau_D^x}^{\tau_D^x}\bigl\langle\nabla u\bigl(X(s\wedge\tau
_D^x)\bigr),dM_x(s)\bigr\rangle=\int_{t\wedge
\tau_D^x}^{\tau_D^x}\langle Z_x(s),dM_x(s)\rangle
\]
for all $t$. By It\^o's isometry,
we have
%
\begin{eqnarray}\label{3.10}
&&E \biggl[ \biggl(\int_{0}^{\infty}\bigl\langle\bigl(\nabla u(X(s
))-Z_x(s)\bigr)\chi_{\{s<\tau_D^x\}},dM_x(s)\bigr\rangle \biggr)^2 \biggr]\nonumber\\
&&\qquad = E \biggl[ \biggl(\int_{0}^{\infty}\bigl\langle\bigl(\nabla u(X(s
))-v_0(X_x(s))\bigr)\chi_{\{s<\tau_D^x\}},dM_x(s)\bigr\rangle \biggr)^2
\biggr]\nonumber\\[-8pt]\\[-8pt]
&&\qquad =E \biggl[ \int_{0}^{\infty}\bigl\langle a(X_x(s))\bigl(\nabla u(X(s
))-v_0(X_x(s))\bigr), \nonumber\\
&&\hspace*{122pt} \nabla u(X(s))-v_0(X_x(s))
\bigr\rangle \chi_{\{s<\tau_D^x\}}\,ds \biggr]=0.\nonumber
\end{eqnarray}
By Fubini theorem and the uniform ellipticity of the matrix $a(x)$,
we deduce that
\[
P^D_s(|\nabla u-v_0|^2)=E\bigl[|\nabla u(X(s))-v_0(X_x(s))|^2
\chi_{\{s<\tau_D^x\}}\bigr]=0
\]
a.e. in $s$ with respect to the Lebesgue
measure, where $P_s^Dh(x)=E_x[h(X(t)),\break t<\tau_D^x]$. The strong
continuity of the semigroup $P^D_s, s\geq0$ implies that
%
\begin{equation}\label{3.11}
|\nabla u-v_0|^2(x)=\lim_{s\rightarrow0}P^D_s(|\nabla
u-v_0|^2)=0
\end{equation}
a.e. Returning to problem (\ref{3.8}), we see that $u$ actually is a
weak solution to the nonlinear problem:
%
\begin{equation}\label{3.12}
\cases{
L_0 u= -f(x, u(x), \nabla u(x)), &\quad in $D$, \cr
u= \phi,&\quad on $\partial D$.
}
\end{equation}
Suppose $\bar{u}$ is another solution to the problem (\ref{3.12}).
By the decomposition of the Dirichlet process $\bar{u}(X_x(s))$, we
find that $(\bar{u}(X_x(s)), \nabla\bar{u}(X_x(s)))$ is also a
solution to the BSDE (\ref{3.7}). The uniqueness of the BSDE
implies that $\bar{u}(X_x(s))=Y_x(s)$. In particular,
$\bar{u}(x)=u_0(x)=Y_x(0)$. This proves the
uniqueness.
\end{pf}

\section{Semilinear elliptic PDEs with singular coefficients}\label{sec5}

In this section, we study the semilinear
second order elliptic PDEs of the following form:
%
\begin{equation}\label{4.1}
\cases{
\mathcal{A}u(x)=-f(x,u(x)), &\quad $\forall x\in D$,\cr
u(x)|_{\partial D}=\varphi, &\quad $\forall x\in\partial{D}$,
}
\end{equation}
where the operator $\mathcal{A}$ is given by
\[
\mathcal{A}={1\over2}\sum_{i,j=1}^{d}{\partial\over{\partial
x_i}}\biggl(a_{ij}(x)\,{\partial\over{\partial x_j}}\biggr)+
\sum_{i=1}^{d}b_i(x)\,{\partial\over{\partial
x_i}}-\mbox{``$\operatorname{div}(\hat{b}\cdot)$''}+q(x)
\]
as in Section~\ref{sec2} and $\varphi\in C(\partial D)$. Consider
the following conditions:
\begin{itemize}[(D.1)]
\item[(D.1)] $(y_1-y_2)(f(x,y_1)-f(x,y_2))\leq-J_1(x)|y_1-y_2|^2$,
\item[(D.2)] $|f(x,y,z)|\leq C,$
\end{itemize}
where $J_1(x)$ is a measurable function, $C$ is a constant.
The following theorem is the main result
of this section.

\begin{theorem}\label{thm5.1}
Suppose that \textup{(D.1)}, \textup{(D.2)} hold and
%
\begin{eqnarray} \label{4.0}
&& E_x^0 \biggl[ \exp \biggl\{\int_0^{\tau_D}\langle (a^{-1}b) (X^0(s)),
dM^0(s)\rangle\nonumber\\
&&\hphantom{E_x^0 \biggl[\exp \biggl\{}
{} + \biggl(\int_0^{\tau_D}\langle (a^{-1}\hat{b}) (X^0(s)), dM^0(s)\rangle
\biggr)\circ
r_{\tau_D}\nonumber\\
&&\hphantom{E_x^0 \biggl[\exp \biggl\{}
{}-\frac{1}{2}\int_0^{\tau_D}(b-\hat{b})a^{-1}(b-\hat{b})^{\ast}(X^0(s))\,ds\\
&&\hphantom{E_x^0 \biggl[\exp \biggl\{}
{}
+\int_0^{\tau_D}q(X^0(s))\,ds-2\int_0^{\tau_D}J_1(X^0(s))\,ds \biggr\} \biggr]\nonumber
\\
&&\qquad <\infty\nonumber
\end{eqnarray}
for some $x\in D$, where $X^0$ is the diffusion generated by $L_0$ as
in Section~\ref{sec2} and $\tau_D$ is the first exit time
of $X^0$ from $D$. Then there exists a unique, continuous weak solution to
equation (\ref{4.1}).
\end{theorem}

\begin{pf}
Set
%
\begin{eqnarray} \label{4.00}
Z_t&=&\exp \biggl\{\int_0^{t}\langle(a^{-1}b) (X^0(s)), dM^0(s)\rangle
+ \biggl(\int_0^{t}\langle(a^{-1}\hat{b}) (X^0(s)), dM^0(s)\rangle
\biggr)\circ
r_{t}\nonumber\\
&&\hphantom{\exp \biggl\{}
{}-\frac{1}{2}\int_0^{t}(b-\hat{b})a^{-1}(b-\hat{b})^{\ast}(X^0(s))\,ds
+\int_0^{t}q(X^0(s))\,ds\\
&&\hspace*{227pt}
{}-2\int_0^{t}J_1(X^0(s))\,ds \biggr\}.\nonumber
\end{eqnarray}
Put
\[
\hat{M}(t)=\int_0^t \langle(a^{-1}\hat{b}) (X^0(s)), dM^0(s)\rangle\qquad
\mbox{for } t\geq0.
\]
Let $R>0$ so that $D\subset B_R:=B(0, R)$. By Lemma 3.2 in \cite{CZ} (see also \cite{Fi}),
there exits a bounded function $v\in W^{1,p}_0(B_R)\subset
W^{1,2}_0(B_R)$ such that
\[
(\hat{M}(t))\circ r_t=-\hat{M}(t)+N^v(t),
\]
where $N^v$ is the zero energy part of the Fukushima decomposition for
the Dirichlet
process $v(X^0(t))$.
Furthermore,
$v$ satisfies the following equation in the
distributional sense:
%
\begin{equation}\label{e:vv}
\operatorname{div}(a \nabla v)=-2 \operatorname{div}(\hat{b}) \qquad \mbox{in } B_R.
\end{equation}
Note that by Sobolev embedding theorem, $v\in C(R^d)$ if we extend
$v=0$ on $D^c$. This implies that $\hat{M}$ and $N^v$ are continuous
additive functionals of $X^0$ in the
strict sense (see \cite{Fi,FOT}), and so is $t\rightarrow(\hat{M}(t))
\circ r_t$.
Thus,
\begin{eqnarray}
&& \biggl(\int_0^{t}\langle(a^{-1}\hat{b}) (X^0(s)), dM^0(s)\rangle
\biggr)\circ
r_{t} \nonumber \\
&&\qquad =-\int_0^{t} \langle(a^{-1}\hat{b}) (X^0(s)), dM^0(s)\rangle+N^v(t)
\nonumber \\
&&\qquad = -\int_0^{t}\langle(a^{-1}\hat{b}) (X^0(s)),
dM^0(s)\rangle+v(X^0(t))-v(X^0(0))-M^v(t)
\nonumber \\
&&\qquad = -\int_0^{t}\langle(a^{-1}\hat{b}) (X^0(s)),
dM^0(s)\rangle+v(X^0(t))-v(X^0(0))\nonumber\\
&&\quad \qquad {} -\int_0^{t}\langle\nabla
v(X^0(s)),dM^0(s)\rangle. \nonumber
\end{eqnarray}
Hence,
\begin{eqnarray}\label{e:Zv}
Z_{t}&=& \frac{e^{v(X^0(t))}}{e^{v(X^0(0))}}\nonumber \\
&&
{}\times\exp \biggl( \int_0^{t} \langle a^{-1}(b -\hat{b}-a\nabla v ) (X^0(s)),
dM^0(s)\rangle-2\int_0^{t}J_1(X^0(s))\,ds \nonumber\\[-8pt]\\[-8pt]
&&\hphantom{{}\times\exp \biggl(}{} +\int_0^{t} \biggl(q-\frac12
(b-\hat{b}-a\nabla v) a^{-1}(b-\hat{b}-a\nabla v)^* \biggr)(X^0(s))\,ds\nonumber
\\
&&\hspace*{95pt}{} + \int_0^{t} \biggl(\frac{1}{2}
(\nabla v)a(\nabla v)^* -\langle b-\hat{b}, \nabla v\rangle \biggr)(X^0(s))\,ds \biggr).\nonumber
\end{eqnarray}
Note that $Z_{t}$ is well defined under $P^0_x$ for every $x\in D$.
Set
$h(x)=e^{v(x)}$. Introduce
\begin{eqnarray*}
\hat{\mathcal{A}}&=&{1\over2}\sum_{i,j=1}^{d}{\partial\over{\partial
x_i}}\,\biggl(a_{ij}(x)\,{\partial\over{\partial x_j}}\biggr)
+\sum_{i=1}^{d}[b_i(x)-\hat{b}_i(x)-(a\nabla
v)_i(x)]\,{\partial\over{\partial x_i}}
\\
&&{}-\langle b-\hat{b},\nabla v\rangle(x)
+{1\over2}(\nabla v)a(\nabla v)^*(x)+q(x).
\end{eqnarray*}
Let $(\Omega, \mathcal{ F}, \mathcal{ F}_t, \hat{X}(t), \hat{P}_x, x\in
R^d)$ be the diffusion process whose infinitesimal generator is given by
\[
\hat{L}={1\over2}\sum_{i,j=1}^{d}{\partial\over{\partial
x_i}}\,\biggl(a_{ij}(x)\,{\partial\over{\partial x_j}}\biggr)
+\sum_{i=1}^{d}[b_i(x)-\hat{b}_i(x)-(a\nabla
v)_i(x)]\,{\partial\over{\partial x_i}}.
\]
It is known from \cite{LZ} that $\hat{P}_x$ is absolutely continuous
with respect to $P_x^0$ and
\[
\frac{d\hat{P}_x}{dP_x^0}\bigg\vert_{\mathcal{ F}_t}=\hat{Z}_t,
\]
where
%
\begin{eqnarray}
\hat{Z}_{t}&=&
\exp \biggl( \int_0^{t}\langle a^{-1}(b -\hat{b}-a\nabla v ) (X^0(s)), dM^0(s)\rangle
\nonumber\\[-8pt]\\[-8pt]
&&\hphantom{\exp \biggl(}{} -\int_0^{t} \biggl(\frac12
(b-\hat{b}-a\nabla v) a^{-1}(b-\hat{b}-a\nabla v)^* \biggr)(X^0(s))\,ds \biggr).\nonumber
\end{eqnarray}
Put
\[
\hat{f}(x,y)=h(x)f(x,h^{-1}(x)y).
\]
Then
\[
(y_1-y_2)\bigl(f(x,y_1)-f(x,y_2)\bigr)\leq-J_1(x)|y_1-y_2|^2.
\]
Consider the following nonlinear elliptic partial differential
equation:
%
\begin{equation}\label{4.2}
\cases{
\hat{\mathcal{A}}\hat{u}(x)=-\hat{f}(x,\hat{u}(x)),&\quad $ \forall x\in
D$,\cr
\hat{u}(x)|_{\partial
D}=h(x)v(x),&\quad $ \forall x\in\partial D$.
}
\end{equation}
In view of (\ref{e:Zv}), condition (\ref{4.0}) implies that
\begin{eqnarray} \quad\nonumber
&& \hat{E}_x \biggl[
\exp \biggl( -2\int_0^{\tau_D}J_1(X^0(s))\,ds +\int_0^{\tau_D}q(X^0(s))\,ds
\nonumber\\[-8pt]\\[-8pt]
&&\hphantom{\hat{E}_x \biggl[\exp \biggl(}
{} + \int_0^{\tau_D} \biggl(\frac{1}{2}
(\nabla v)a(\nabla v)^* -\langle b-\hat{b}, \nabla v\rangle \biggr)(X^0(s))\,ds\biggr) \biggr]<\infty,\nonumber
\end{eqnarray}
where $\hat{E}_x$ indicates that the expectation is taken under $\hat{P}_x$.
From Theorem~\ref{thm4.2}, it follows that equation (\ref{4.2}) admits a unique
weak solution $\hat{u}$. Set $u(x)=h^{-1}(x)\hat{u}(x)$. We will
verify that $u$ is a weak solution to equation (\ref{4.1}).

Indeed, for $\psi\in W^{1,2}_{0}(D)$, since $\hat{u}(x)=h(x)u(x)$ is a
weak solution to equation
(\ref{4.2}), it follows that
\begin{eqnarray*}
&&{1\over2}\sum_{i,j=1}^{d}\int_Da_{ij}(x)\,{\partial
[h(x)u(x)]\over{\partial x_i}}\,{\partial[h^{-1}(x)\psi]\over
{\partial x_j}}\,dx\\
&&\quad
{}-\sum_{i=1}^{d}\int_D[b_i(x)-\hat{b}_i(x)-(a\nabla
v)_i(x)]\,{\partial[h(x)u(x)]\over{\partial x_i}}\,h^{-1}(x)\psi
\,dx\\
&&\quad {}+\int_D\langle b-\hat{b},\nabla v(x)\rangle u(x)\psi(x)\,dx \\
&&\quad {}-{1\over2}\int_D(\nabla v)a(\nabla v)^*(x)u(x)\psi \,dx-\int
_Dq(x)u(x)\psi(x) \,dx\\
&&\qquad =\int_Df(x, u(x))\psi(x) \,dx.
\end{eqnarray*}
Denote the terms on the left of the above equality, respectively, by
$T_1$, $T_2$, $T_3$, $T_4$, $T_5$. Clearly,
%
\begin{eqnarray}\label{4.3}
T_1&=& {1\over2}\sum_{i,j=1}^{d}\int_Da_{ij}(x)\,{\partial
u(x)\over{\partial x_i}}\,{\partial\psi\over
{\partial x_j}}\,dx-{1\over2}\sum_{i,j=1}^{d}\int_Da_{ij}(x)\,{\partial
u(x)\over{\partial x_i}}\,{\partial v\over
{\partial x_j}}\,\psi \,dx\nonumber\\
&&{}+ {1\over2}\sum_{i,j=1}^{d}\int_Da_{ij}(x)\,{\partial
v\over{\partial x_i}}\,{\partial\psi\over
{\partial x_j}}\,u(x)\,dx\\
&&{}-{1\over2}\sum_{i,j=1}^{d}\int
_Da_{ij}(x)\,{\partial
v\over{\partial x_i}}\,{\partial v\over
{\partial x_j}}\,\psi u(x)\,dx.\nonumber
\end{eqnarray}
Using chain rules, rearranging terms, it turns out that
%
\begin{eqnarray}\label{4.4}
T_2+T_3&=&-\sum_{i=1}^{d}\int_Db_i(x)\,{\partial u(x)\over\partial
x_i}\,\psi\, dx-\sum_{i=1}^{d}\int_D\hat{b}_i(x){\partial\psi\over{\partial
x_i}}\,u(x)\,dx\nonumber\\
&&{}+\sum_{i=1}^{d}\int_D[\hat{b}_i(x)+(a\nabla
v)_i(x)]\,{\partial[\psi u(x)]\over{\partial x_i}}\,dx\\
&&{}-\sum_{i=1}^{d}\int_D(a\nabla
v)_i(x)\,{\partial\psi\over{\partial x_i}}\,u(x) \,dx+\sum_{i=1}^{d}\int
_D(a\nabla
v)_i(x)\,{\partial v \over{\partial x_i}}\,u(x)\psi \,dx.\nonumber
\end{eqnarray}
In view of (\ref{e:vv}),
%
\begin{eqnarray} \label{4.5}
&&\sum_{i=1}^{d}\int_D[\hat{b}_i(x)+(a\nabla
v)_i(x)]\,{\partial[\psi u(x)]\over{\partial x_i}}\,dx\nonumber\\[-8pt]\\[-8pt]
&&\qquad =\frac{1}{2}\sum
_{i=1}^{d}\int_D(a\nabla
v)_i(x)\,{\partial[\psi u(x)]\over{\partial x_i}}\,dx.\nonumber
\end{eqnarray}
Thus,
%
\begin{eqnarray}
T_2+T_3
&=&-\sum_{i=1}^{d}\int_Db_i(x)\,{\partial u(x)\over\partial x_i}\,\psi
\,dx-\sum_{i=1}^{d}\int_D\hat{b}_i(x)\,{\partial\psi\over{\partial
x_i}}\,u(x)\,dx\nonumber\\
&&{}+\frac{1}{2}\sum_{i=1}^{d}\int_D(a\nabla
v)_i(x)\,{\partial[\psi u(x)]\over{\partial x_i}}\,dx\\
&&{}-\sum_{i=1}^{d}\int_D(a\nabla
v)_i(x)\,{\partial\psi\over{\partial x_i}}\,u(x) \,dx+\sum_{i=1}^{d}\int
_D(a\nabla
v)_i(x)\,{\partial v \over{\partial x_i}}\,u(x)\psi \,dx.\nonumber
\end{eqnarray}

After cancelations, it is now easy to see that
%
\begin{eqnarray}
T_1+T_2+T_3+T_4+T_5
&=&{1\over2}\sum_{i,j=1}^{d}\int_Da_{ij}(x)\,{\partial u(x)\over
{\partial x_i}}\,{\partial\psi\over{\partial
x_j}}\,dx \nonumber\\
&&{}-\sum_{i=1}^{d}\int_Db_i(x)\,{\partial u(x)\over
{\partial x_i}}\,\psi \,dx\nonumber\\[-8pt]\\[-8pt]
&&{}-\sum_{i}^{d}\int_D\hat{b}\,{\partial
\psi\over
\partial x_i}\,u(x)\,dx-\int_Dq(x)u(x)\psi(x) \,dx\nonumber \\
&=&\int_Df(x, u(x))\psi(x) \,dx.\nonumber
\end{eqnarray}
Since $\psi$ is arbitrary, we conclude that $u$ is a weak solution of
equation (\ref{4.1}). Suppose $u$ is a continuous weak solution to
equation (\ref{4.1}). Put $\hat{u}(x)=h(x)u(x)$. Reversing the above
process, we see that
$\hat{u}$ is a weak solution to equation (\ref{4.2}). The uniqueness of
the solution of equation (\ref{4.1}) follows from that of equation (\ref{4.2}).
\end{pf}

%

\printaddresses

\end{document}